\documentclass[smallextended]{svjour3}       

\usepackage{graphicx,float}
\usepackage{mathptmx}      
\DeclareMathAlphabet{\mathcal}{OMS}{cmsy}{m}{n}  
\usepackage{amsmath}
\usepackage{amssymb}
\usepackage{appendix}
\usepackage{epsfig}
\usepackage{hyperref}
\usepackage{stackengine}
\stackMath

\makeatletter

\renewcommand*\env@matrix[1][\arraystretch]{%
  \edef\arraystretch{#1}%
  \hskip -\arraycolsep
  \let\@ifnextchar\new@ifnextchar
  \array{*\c@MaxMatrixCols c}}

\newcommand{\stabfigurepath}{.}
\newcommand{\plotfigurepath}{.}

\journalname{Noname}

\begin{document}

\title{A parallel-in-time approach for wave-type PDEs}

\author{Abe C. Ellison \and Bengt Fornberg}

\institute{
    Department of Applied Mathematics \\
    University of Colorado \\
    Boulder, CO 80309 \\
    \email{abe.ellison@colorado.edu}
}

\date{Received: date / Accepted: date}

\maketitle

\begin{abstract}
Numerical solutions to wave-type PDEs utilizing method-of-lines require
the ODE solver's stability domain to include a large stretch of the
imaginary axis surrounding the origin.  We show here that extrapolation 
based solvers of Gragg-Bulirsch-Stoer (GBS) type can meet this requirement.
Extrapolation methods utilize several independent time stepping sequences,
making them highly suited for parallel execution.  Traditional extrapolation schemes 
use all time stepping sequences to maximize the method's order of accuracy.
The present method instead maintains a desired order of accuracy while employing
additional time stepping sequences to shape the resulting stability domain.
We optimize the extrapolation coefficients to maximize the stability domain's imaginary axis coverage.
This yields a family of explicit schemes that approaches maximal time step size for wave propagation problems.
On a computer with several cores we achieve both high order and fast time to solution compared
with traditional ODE integrators.
\end{abstract}

\section{Introduction}

Time integration of ODEs is an inherently sequential process, since
each forward step ought to be based on the most recent information
available. Three conceivable options for achieving some level of \emph{parallel-in-time}
are (i) to have correction calculations follow the explicit forward
steps as closely behind as possible, letting them catch up frequently,
(ii) to carry out `preparatory' calculations that are based on trying
to anticipate later solution states, and (iii) to exploit extrapolation
ideas. While all of these concepts have been pursued for systems of
ODEs, as summarized in \cite{CMO10,KW14}, their performance is unclear
for ODE systems arisen from method-of-lines (MOL) discretization of
wave-type PDEs. The additional requirement that arises
then is that the ODE solver's \emph{stability domain} must include
a quite large stretch of the imaginary axis surrounding the origin.
We show here that extrapolation-based ODE solvers of Gragg-Bulirsch-Stoer
(GBS) type can meet this requirement. In particular, one such scheme
that we will focus on steps forward explicitly using six cores as
fast as Forward Euler (FE) does on one core, but combines eighth order
of accuracy with a generously sized stability domain. In contrast
to linear multistep methods, it needs no back levels in time to get
started. The present approach is compared against explicit Runge-Kutta (RK) methods
for a PDE test problem. 

Standard Richardson extrapolation schemes utilize a square Vandermonde-type system
to compute the extrapolation weights.  This system is constructed
to cancel successive terms in the asymptotic error expansion of the time stepper.
By allowing more columns than rows in the system - that is, more extrapolation
components than order constraints - we create an underdetermined system
that grants degrees of freedom to optimize the extrapolated stability domain.
For wave-type PDEs stepped with method-of-lines we optimize the stability
domain along the imaginary axis.  We achieve stability domains far larger
than those of both the square extrapolation systems and other standard ODE integrators,
thereby enabling large time step sizes and thus faster time to solution.

\section{GBS-type ODE solvers}

\subsection{GBS concept}

We consider first the problem of advancing forward in time an ODE of
the form $y'(t)=f(t,y)$, where the unknown function $y(t)$ is either
scalar or vector valued. The complete time interval of interest is
split into $N$ sections. For each of these sections, the
basic (unextrapolated) GBS scheme consists of the steps
\begin{equation}
\left\{
    \bgroup
    \renewcommand*{\arraystretch}{2.2}
    \begin{array}{ll}
        {\displaystyle \frac{y_{1}-y_{0}}{h}}=f(t_{0},y_{0}) & \textrm{Forward Euler (FE)}\\
        {\displaystyle \frac{y_{n+1}-y_{n-1}}{2h_{\phantom{}}}}=f(t_{n},y_{n}) & \textrm{Leap-frog (LF), }n=1,2,\ldots,N\\
        y_{N}^{*}=\dfrac{1}{4}(y_{N-1}+2y_{N}+y_{N+1}) & \textrm{Averaging},
    \end{array}
    \egroup
\right.
\label{eq:GBS}
\end{equation}
after which $y_{N}^{*}$ is accepted as the new value at time $t_{N}.$
The initial FE step is accurate to first order, while the subsequent
LF steps are second order accurate. One would therefore expect $y_{N}^{*}$
to be accurate to at most second order, and have an error expansion
in which all further powers of $h$ would be present. Remarkably,
for any smooth (linear or nonlinear) function $f(t,y)$, it transpires
that, if $N$ is even, all odd powers in the expansion will vanish
\cite{BS66,G65,HNW87,L91}:
\begin{align}
\text{Error}=y_{N}^{*}-y(t_{N})=a_{2}h^{2}+a_{4}h^{4}+a_{6}h^{6}+\ldots
\end{align}
The form of the expansion makes Richardson extrapolation
particularly efficient since, each time this is applied, the result
will gain two orders of accuracy. For example, the results from four
completely independent calculations over the same section in time,
using different $h$-values, can be combined to give an $\mathcal{O}(h^{8})$-accurate
result. These four calculations require no communications between
each other, and can therefore be run simultaneously on separate cores.

Even when not counting the cost of the work on the extra cores, the GBS approach does not
offer any striking benefits for standard ODE systems, unless possibly
if extrapolated to very high orders. However, in the present context
of wave-type PDEs, the situation becomes different, since GBS methods
can be designed to feature particularly favorable stability domains. 

\subsection{Stability domains for GBS-type methods}

\ref{sec:appendix_a} briefly summarizes the definition of an ODE solver's
stability domain, explains its significance in the context of
MOL time stepping, and provides stability domain information for some
well-known explicit ODE solvers. These domains should be contrasted
to the corresponding ones for GBS methods described below. The imaginary
stability boundary (ISB) of an ODE solver is defined as the largest
value such that the imaginary axis is included from $-i\cdot \text{ISB}$
to $+i\cdot \text{ISB}$. For solvers that lack any imaginary axis coverage,
we define their ISB to be zero. 

In order to provide a fair comparison between different methods, we
will from now on further normalize all ISB values by the number of
function evaluations that each step requires and denote this ISB$_{\text{n}}$.
For example, we divide
$\textrm{R}\textrm{K}_{4}$'s ISB, stated in \ref{sec:appendix_a} as 2.8284,
by four to compensate for its four stages (with one function evaluation
in each), i.e. we list its ISB$_{\text{n}}$
as 0.7071 $(=1/\sqrt{2}$). Similarly, the ISB$_{\text{n}}$
for the 13-stage $\textrm{R}\textrm{K}_{8}$ method becomes 0.2848.
With this normalization, the largest feasible ISB$_{\text{n}}$
for any explicit method becomes one \cite{JN81}, which is realized
by the LF scheme. Since the longest distance a solution can be advanced
per function evaluation is proportional to the time stepping method's
ISB$_{\text{n}}$, a key goal will be to design a method that has both
high order and a large ISB$_{\text{n}}$. 

Stability domains and ISBs for GBS-type schemes do not appear to have been
studied until in \cite{FZL07}. One key observation made there was that GBS schemes of
orders 4, 8, 12, ... will feature positive ISBs, whereas schemes of order 6,
10, 14, ... will not.  Hence, in what follows we study only schemes with order
divisible by four.

\section{Optimizing the Stability Domain}

\subsection{Introduction to ISB Optimization}

Stability domain optimization has been well studied in the literature. 
The class of steppers that maximizes ISB given stability polynomial order $N+1$ was found
independently by Kinnmark and Gray \cite{KG84a} and  Sonneveld and van Leer \cite{SL85}.
The methods divide a time interval into $N$ evenly spaced steps.  A Forward Euler predictor
and Backward Euler corrector pair initiates the time step, then $N-1$ leap frogs
bring us to the end of the time interval.  This class of methods has order of accuracy at most
two, and achieves an ISB$_{\text{n}}=N/(N+1)$.

Kinnmark and Gray demonstrate third
and fourth order accurate stability polynomials in \cite{KG84b} that converge
to the optimal ISB as number of subintervals increases.  Interestingly, the first
two methods of this class are the third order and fourth order explicit Runge-Kutta methods
with three and four stages, respectively.  Thus RK$_4$ is optimal in the sense that it
fully utilizes its four function evaluations to maximize time step for wave-type problems.
It is therefore an excellent candidate for comparison with the optimized schemes that follow.

\subsection{GBS Stability Domain Optimization}

In Richardson extrapolation schemes one sets up a square Vandermonde
system to compute the weights guaranteeing a specified order of accuracy.
If we allow the number of components in the extrapolation scheme to
increase beyond those necessary for maintaining order of accuracy
we obtain an underdetermined system with degrees of freedom. We
utilize these degrees of freedom to optimize the stability domain
along a contour in the complex plane. 

Extrapolation allows us to eliminate successively higher order terms
in the asymptotic error expansion of our solution. To do so, for each of \emph{$m$}
integrators we divide the time interval $H$ into $n_{i}$ steps of size $h_i=H/n_i$, $i=1,2,\hdots,m$.  We then
construct a linear system to eliminate terms through order $p-1$ in the error expansion,
yielding a $p$-order accurate solution.
In the case of GBS integrators, the odd coefficients in the asymptotic expansion
are zero.  Thus we may drop the constraint equations for odd powers
of $h_i$, obtaining a system of $\frac{p}{2}$ equations:
\begin{align} 
\underbrace{\addstackgap[16pt]{
  \begin{bmatrix}
    1 & 1 & \hdots & 1 \\
    h_1^2 & h_2^2 & \hdots & h_m^2 \\
    \vdots & \vdots & \ddots & \vdots \\
    h_1^{p-2} &  h_2^{p-2} & \hdots & h_m^{p-2}
  \end{bmatrix}_{\frac{p}{2}\times m}
}}_{V}
\hspace{.5ex}
\underbrace{\addstackgap[5pt]{
  \begin{bmatrix}[1]
    c_1 \\ c_2 \\ \\ \vdots \\ \\ c_m
  \end{bmatrix}_m
}}_{c}
=
\underbrace{\addstackgap[16pt]{
  \begin{bmatrix}[1]
    1 \\
    0 \\
    \vdots \\
    0
  \end{bmatrix}_{\frac{p}{2}}.
}}_{b}
\label{eq:vander}
\end{align}
When $m=\frac{p}{2}$ we have a square matrix which corresponds to
the usual Richardson extrapolation schemes. The matrix is invertible
when $n_{i}\ne n_{j},i\ne j$, and so we solve for the weight vector $c$,
which we apply to the individual integrated solutions to form the combined
solution at the end of the time interval.

By allowing $m>\frac{p}{2}$ the system becomes underdetermined and we may
enforce order constraints while optimizing selected features of the stability
domain.  The optimization algorithm was adapted from the polynomial optimization formulation
in \cite{KA12}; details are provided in \ref{sec:appendix_b}.

\subsection{Fully-Determined Optimization Results}

We first investigate optimal step count selection for fully determined
extrapolation schemes.  For schemes of order $p$ we test each combination of $\frac{p}{2}$ step counts
up to a set maximum, here chosen to be 24.  Each combination yields a set of uniquely determined 
extrapolation weights.  We then select the combination of step counts
that maximizes ISB$_{\text{n}}$ of the extrapolated stability domain.  Tab.~\ref{tab:isb_table_fulldet} contains the tabulated
results for orders eight, twelve and sixteen.  The schemes all have generous imaginary
axis coverage and can be implemented efficiently on three, four and five cores respectively.
\begin{table}[H]
\caption{Optimal step count sequences and ISB$_{\text{n}}$ for the fully-determined schemes}
\begin{tabular}{|c|c|c|c|c|c|}
\hline
Order & Cores & Step Counts & ISB$_{\text{n}}$ \\
\hline
8  & 3 & 2,16,18,20 & 0.5799 \\
12 & 4 & 2,8,12,14,16,20 & 0.4515 \\
16 & 5 & 2,8,10,12,14,16,18,22 & 0.4162 \\
\hline
\end{tabular}
\label{tab:isb_table_fulldet}
\end{table}

\subsubsection{Eighth Order}\label{sec:discrete_8}
The eighth order, three-core scheme has ISB$_{\text{n}} = 0.5799$ with the following step counts and 
uniquely determined weights:
\begin{equation*}
\begin{aligned}
\begin{matrix}[1.5] \text{Step Counts } \\ \text{Weights } \end{matrix} & 
\begin{matrix}[1.5] : \\ : \end{matrix}
\hspace{1ex}
\begin{matrix}[1.5]
    2, & 16, & 18, & 20 \\
    -\frac{1}{498960}, & \frac{65536}{9639}, & -\frac{531441}{25840}, & \frac{250000}{16929}.
\end{matrix}
\end{aligned}
\end{equation*}

\subsubsection{Twelfth Order}
The twelfth order, four-core scheme has ISB$_{\text{n}} = 0.4515$ with the following step counts and weights:
\begin{equation*}
\begin{aligned}
\begin{matrix}[1.5] \text{Step Counts } \\ \text{Weights } \end{matrix} & 
\begin{matrix}[1.5] : \\ : \end{matrix}
\hspace{1ex}
\begin{matrix}[1.5]
    2, & 8, & 12, & 14, & 16, & 20 \\
    -\frac{1}{157172400},  & \frac{4096}{155925},  & -\frac{59049}{15925},  & \frac{282475249}{15752880},  & -\frac{4194304}{178605},  & \frac{9765625}{954261}.
\end{matrix}
\end{aligned}
\end{equation*}

\subsubsection{Sixteenth Order}
The sixteenth order, five-core scheme has ISB$_{\text{n}}=0.4162$ and utilizes the step count
sequence $\{2,8,10,12,14,16,18,22\}$.  Extrapolation weights can be computed by solving
the corresponding Vandermonde system (\ref{eq:vander}).  We omit them here since the weights are ratios of
large integers in both the numerators and denominators.

\subsection{Underdetermined Optimization Results}\label{sec:optresults}

\begin{figure}[H]
\includegraphics[width=\linewidth]{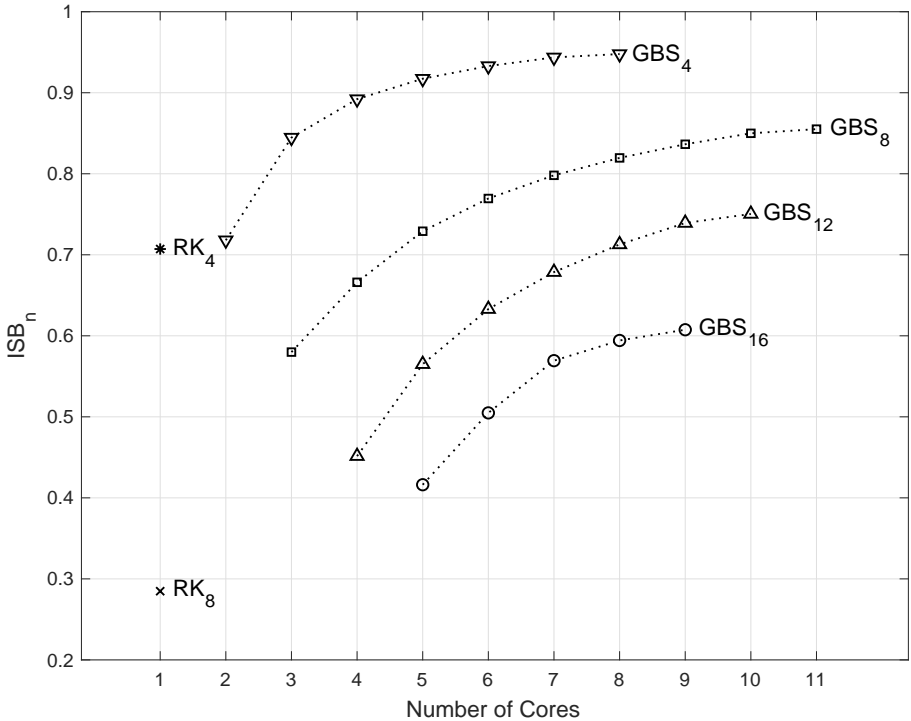}
\caption{Trends for optimized ISB$_{\text{n}}$ versus core count for methods of various order}
\label{fig:isb_vs_cores}
\end{figure}

Using the optimization methodology described in \ref{sec:appendix_b} we optimize the ISB of GBS-type methods up to order
sixteen.  Increasing the number of extrapolation components leads to an
increase in ISB$_{\text{n}}$.  Since for explicit schemes
the maximum ISB$_{\text{n}}$ is one we expect the relative gains of adding more components to eventually saturate.
By evenly distributing work across CPU cores we can demonstrate the relationship
between available processors and maximal time step size.
This correspondence between core count and ISB$_{\text{n}}$ is shown in Fig.~\ref{fig:isb_vs_cores}.
Here we observe that efficiency saturation occurs around ten cores for all orders of accuracy;
the saturation value itself is strongly dependent on the order of accuracy.

\begin{figure}[H]
\parbox{.495\linewidth}{
\includegraphics[width=\linewidth]{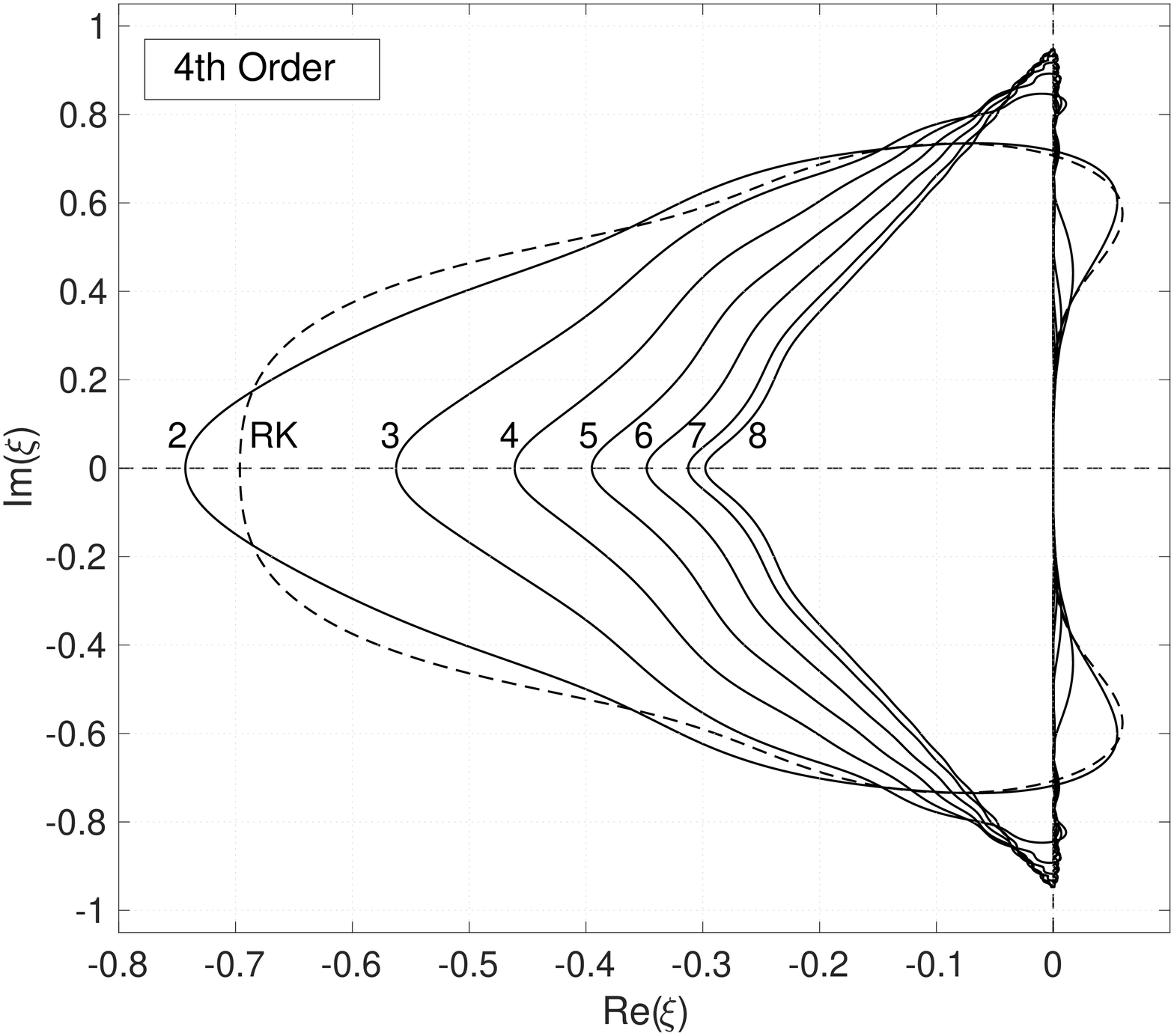}
}
\parbox{.495\linewidth}{
\includegraphics[width=\linewidth]{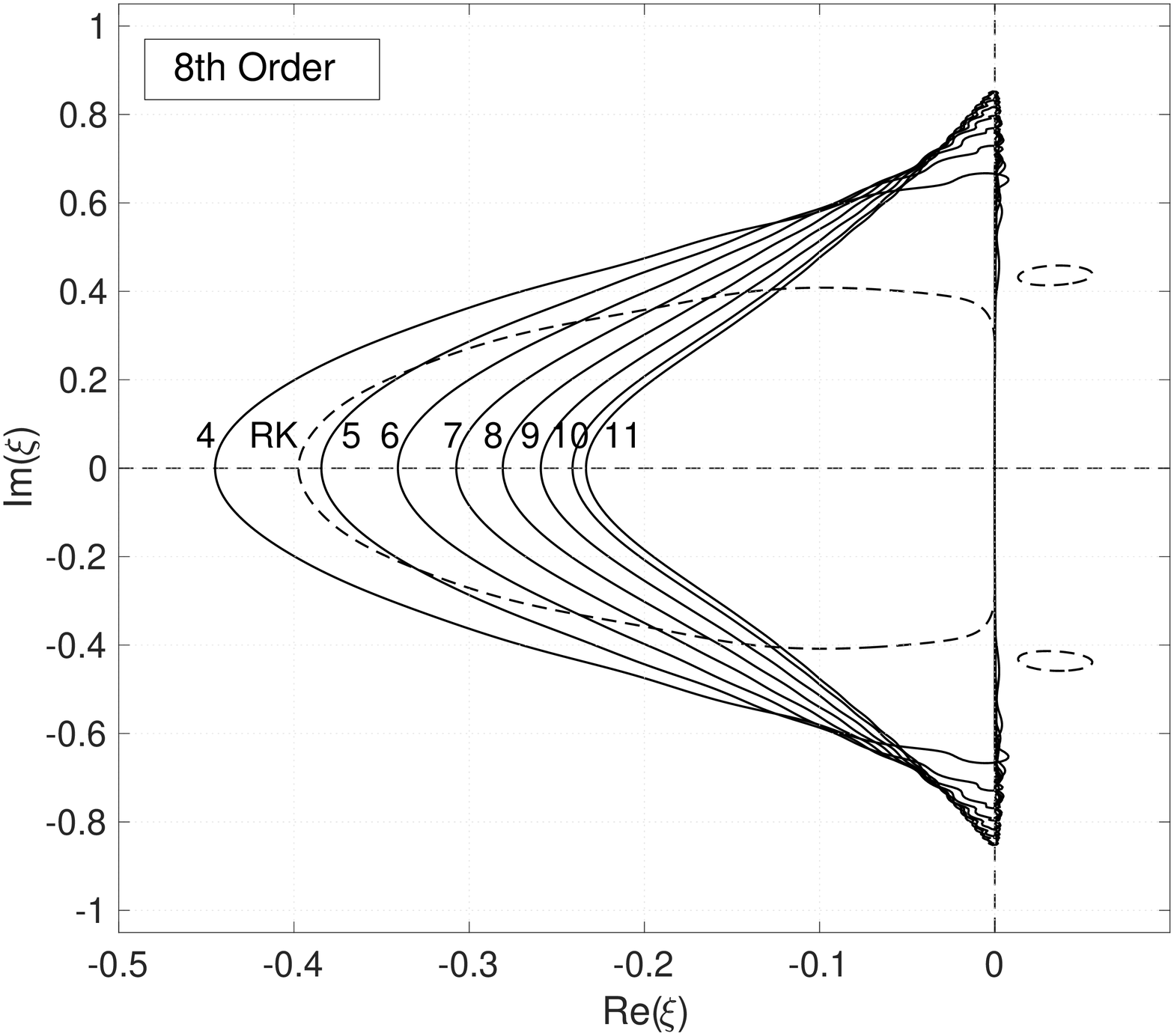}
}\\
\parbox{.495\linewidth}{
\includegraphics[width=\linewidth]{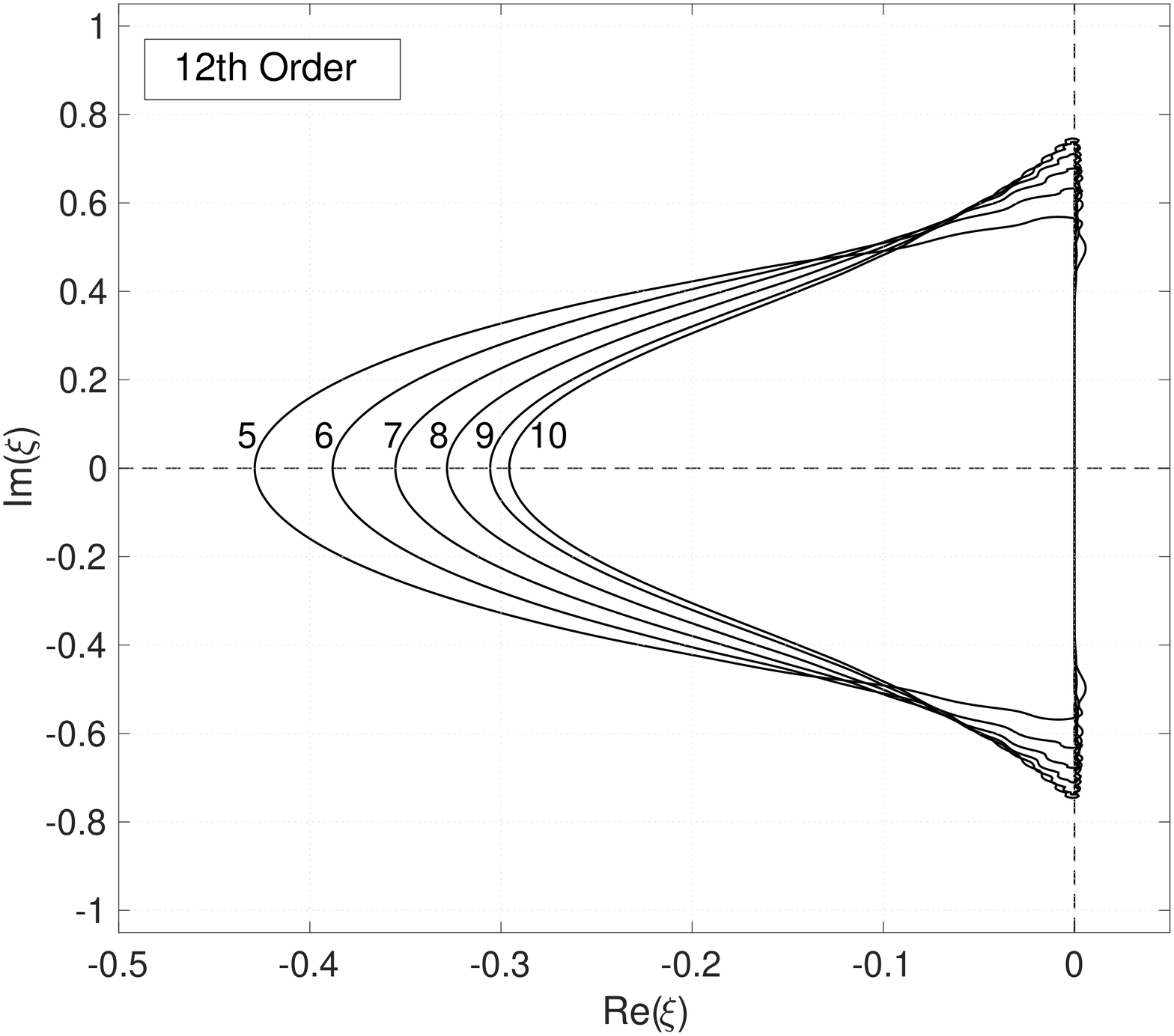}
}
\parbox{.495\linewidth}{
\includegraphics[width=\linewidth]{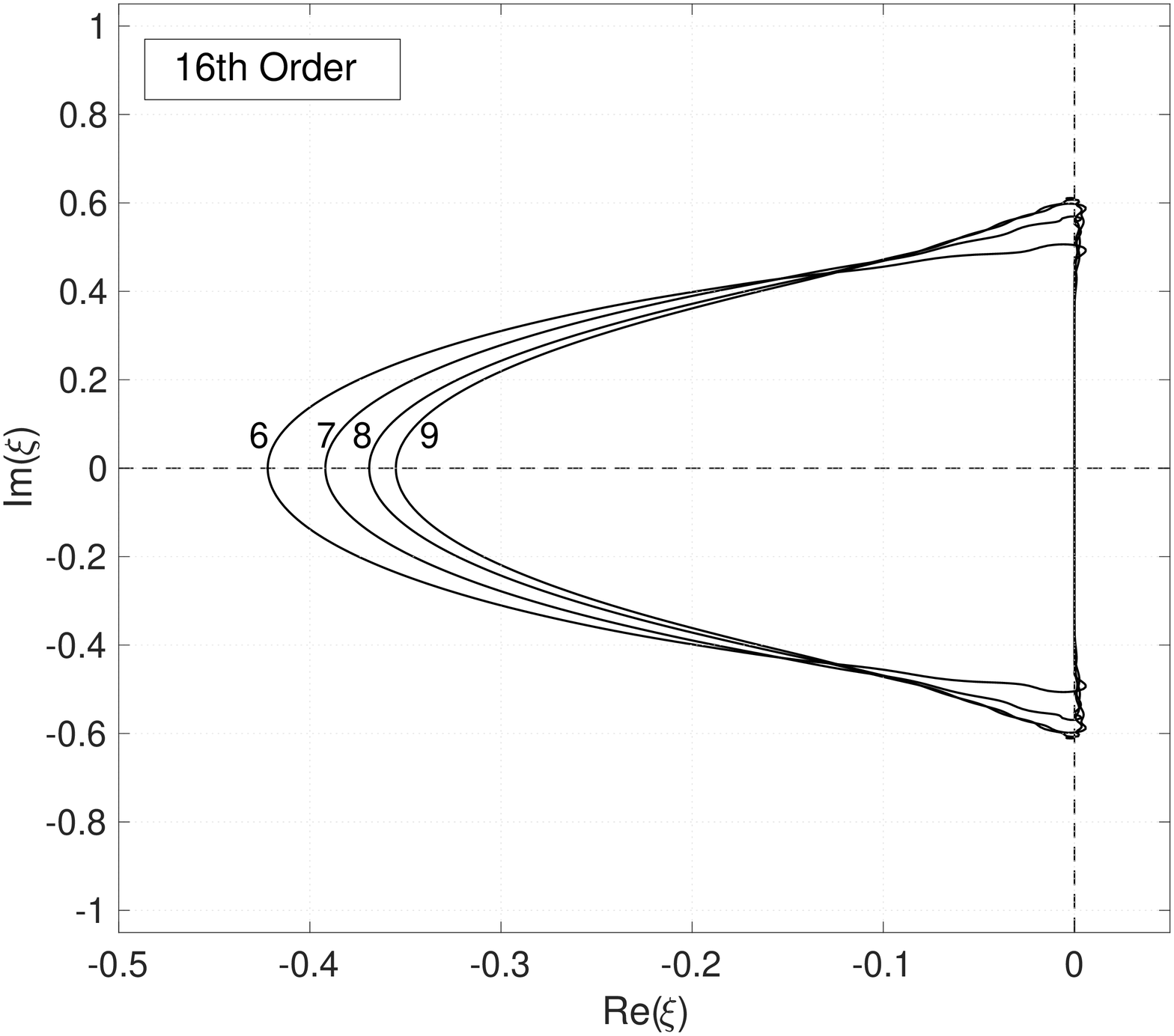}
}
\caption{Optimized stability domains labeled with number of cores required to implement the schemes,
with fourth and eighth order Runge-Kutta (RK) stability domains for comparison}
\label{fig:collected_stabdomains}
\end{figure}

We achieve the best optimization results by utilizing all (even) step counts up to a maximum
dependent on the number of available CPU cores.  The $m$ extrapolation
components then have subinterval counts $\{2, 4, \hdots, 2m\}$, with $m$ set by the available processing resources.
We denote the number of subintervals at which ISB$_{\text{n}}$ can no longer be increased $N_{\text{opt}}$ and
aggregate the results for each order of accuracy in Tab.~\ref{tab:isb_table}.
The stability domains for each core count are plotted
in Fig.~\ref{fig:collected_stabdomains}.
Results demonstrate a tradeoff between optimal ISB and order of accuracy as is typical
of explicit time integrators.
\begin{table}[H]
\caption{Largest optimized ISB$_{\text{n}}$ obtained for the underdetermined schemes}
\begin{tabular}{|c|c|c|c|}
\hline
Order & Cores & $N_{\text{opt}}$ & ISB$_{\text{n}}$ \\
\hline
4  &  8 & 28 & 0.9477 \\
8  & 11 & 40 & 0.8551 \\
12 & 10 & 36 & 0.7504 \\
16 &  9 & 32 & 0.6075 \\
\hline
\end{tabular}
\label{tab:isb_table}
\end{table}

The capping of $N_{\text{opt}}$ is an artifact of the optimizer.  Our convex solver
fails to produce methods with larger ISB if we increase the number of free variables
beyond those presented in Tab.~\ref{tab:isb_table}.  We believe that by addressing the conditioning
of the Vandermonde system as in \cite{KA12} one can continue further along the curves presented in
Fig.~\ref{fig:isb_vs_cores}.  Extrapolating these curves shows the schemes do not converge
to the optimal ISB$_{\text{n}}=1$; the exact tradeoff between order of accuracy and optimal
ISB is a topic of future research.

\subsection{Leading Order Error}
\begin{figure}[H]
\includegraphics[width=\linewidth]{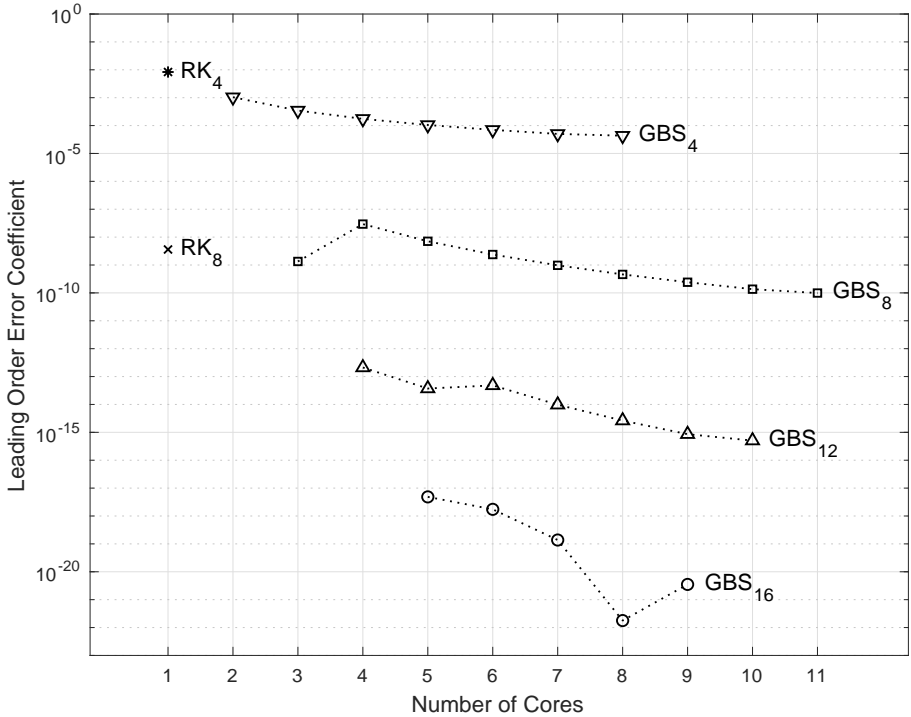}
\caption{Leading order error term $a_{p+1}$ for the optimized methods}
\label{fig:leading_order_error}
\end{figure}

Let the time integrator's stability polynomial, as defined in \ref{sec:appendix_a}, be denoted $R(\xi)$.
For a method of order $p$, the stability domain's boundary follows the imaginary axis
surrounding the origin linearly up to deviation on the order $p+1$.  To compute
the leading error coefficient we set the stability polynomial $R(\xi) = e^{i\theta}$.  Taking the
complex logarithm of the polynomial and Taylor expanding yields a power series for $\theta(\xi)$.
We then compute the inverse series to find
$\xi(\theta) = i\theta+a_{p+1}(i\theta)^{p+1}+a_{p+2}(i\theta)^{p+2}+\mathcal{O}\left(\theta^{p+3}\right)$.
Since we consider only methods with order divisible by four we simplify as follows:
\begin{align}
\xi(\theta) = i\theta + i a_{p+1} \theta^{p+1} - a_{p+2} \theta^{p+2} + \mathcal{O}\left(\theta^{p+3}\right).
\end{align}
Departure from the imaginary axis is governed by the $a_{p+2}$ coefficient.  We then require $a_{p+2}$ to be
negative to ensure the stability domain has a positive ISB.  Accuracy is determined by
the $a_{p+1}$ coefficient; Fig.~\ref{fig:leading_order_error} presents this coefficient
for each method as a function of number of cores.
The leading order error coefficient decays toward zero as number of cores
increases while holding order fixed.  This implies the underdetermined extrapolation
schemes do not trade away numerical precision for achieving large ISBs - they instead
gain significantly in accuracy.

\subsection{Core Partitioning}\label{sec:corepartition}
In order to achieve the theoretical efficiencies presented in Sect.~\ref{sec:optresults} we require 
a work partitioning scheme that distributes the individual time steppers amongst the cores.
For a specified maximum subinterval count $N_{\text{max}}$ we achieve the largest ISB$_{\text{n}}$ by utilizing all even step counts up
to $N_{\text{max}}$.  For a fixed number of cores, here denoted $n_{\text{cores}}$, we can always evenly distribute
the work when we choose $N_{\text{max}}=4\cdot n_{\text{cores}}-2$.  This corresponds to folding together onto a single core
pairs of integrators with step counts adding to $N_{\text{max}}$.  For example, with $N_{\text{max}}=10$,
we evenly load three cores with step counts \{10, \{8, 2\}, \{6, 4\}\}.

The GBS scheme requires $N+1$ function evaluations for an integrator with $N$ subintervals.
When stacking multiple integrators on a single core we share the first function evaluation
since it takes identical arguments for \emph{all} time steppers.  We can in principle share this evaluation
among all cores but communication overhead may make this approach less efficient.

\subsection{Method Specification}\label{sec:notation}
Let the $m\times\frac{p}{2}$-sized Vandermonde matrix in (\ref{eq:vander}) be denoted $V$.
When we have $m>\frac{p}{2}$ the system is underdetermined.  We
then choose $\frac{p}{2}$ dependent step counts to meet the order constraints
and collect these into the step count sequence $\{n_{\text{dep}}\}$.  The associated extrapolation
weight vector is denoted $c_{\text{dep}}$ and the corresponding columns of $V$ denoted $V_{\text{dep}}$.
Likewise collect the remaining step counts into the $\left(m-\frac{p}{2}\right)$-length step count sequence $\{n_{\text{free}}\}$
and label corresponding extrapolation weights $c_{\text{free}}$ and Vandermonde columns $V_{\text{free}}$.
Then the constraint equations may be written
\begin{align}
V_{\text{dep}}c_{\text{dep}}+V_{\text{free}}c_{\text{free}} = b.
\label{eq:order_constraints}
\end{align}
The Vandermonde system order constraints must be met exactly.
We thus omit floating-point coefficients for $c_{\text{dep}}$ in the text; they are best computed
symbolically then converted to the desired floating point format.  We may readily compute $c_{\text{dep}}$
with the relation 
\begin{align}
c_{\text{dep}} = V_{\text{dep}}^{-1}\left(b-V_{\text{free}}c_{\text{free}}\right),
\label{eq:cdep_computation}
\end{align}
so one only needs the step count sequences $\{n_{\text{dep}}\}$ and $\{n_{\text{free}}\}$ and  weight vector
$c_{\text{free}}$ to fully specify a scheme.

\subsection{Methods of Choice}
In this section we present two eighth order methods and one twelfth order method
with rational coefficients for convenient use.  In order to provide
robustness to spurious discretized eigenvalues sitting slightly in
the right-half plane we push the optimization curve into the positive
reals.  This disturbs the ISB$_{\text{n}}$ very little and makes the method 
suitable for local differentiation stencils generated for example by
RBF-FD (radial basis function-generated finite difference) approximations
\cite{FFBook}.  

To generate the following schemes we first optimize the free coefficients using
the optimization methodology described in \ref{sec:appendix_b}.  The optimization contour is
chosen to trade off a small amount of imaginary axis coverage for an area containing
the positive reals away from the origin.  We then perform a search
over a set of rational numbers that closely approximate the floating point free
coefficients.  We select a set with small integers in the numerator and denominator
which disturbs the scheme's stability domain very little.  These coefficients are
reported below.

\subsubsection{The Eighth Order, Six Core Method: GBS$_{8,6}$}
\begin{figure}[H]
\includegraphics[width=\linewidth]{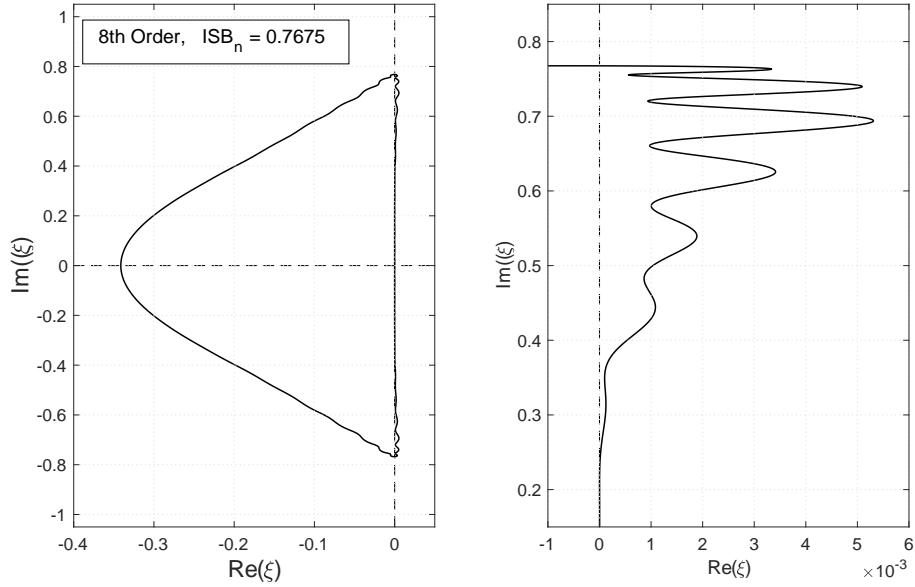}
\caption{Stability domain of the eighth order, six core scheme (left), and zoom (right)}
\label{fig:stabdomain_p8_6core}
\end{figure}

The following eighth order method achieves a robust stability domain for wave-type PDEs
on six cores and is therefore dubbed GBS$_{8,6}$.
The scheme's ISB$_{\text{n}}$ is 0.7675, a 0.26\% reduction from the optimal six core value of 0.7695.
The scheme can be implemented using the following step counts and extrapolation weights:
\begin{equation*}
\begin{aligned}
\{n_{\text{dep}}\} &= \hspace{.3ex} \left\{\begin{matrix} 2, & 4, & 6, & 10 \end{matrix}\right\} \\
\begin{matrix}[1.5] \{n_{\text{free}}\} \\ \vphantom{\bigg[}c_{\text{free}} \end{matrix} & 
\begin{matrix}[1.5] = \\ \vphantom{\bigg[}= \end{matrix}
\hspace{1ex}
\begin{matrix}[1.5]
    \{ & 8, & 12, & 14, & 16, & 18, & 20, & 22 &\}\hphantom{^T.} \\
    \bigg[ & \frac{2165}{767488}, & \frac{13805}{611712}, & \frac{4553}{72080}, & \frac{14503}{66520}, & \frac{27058}{7627}, & -\frac{86504}{5761}, & \frac{40916}{3367} &\bigg]^T.
\end{matrix}
\end{aligned}
\end{equation*}
The dependent weights can be computed
exactly by inverting the corresponding Vandermonde system.
The scheme's stability domain is plotted in Fig.~\ref{fig:stabdomain_p8_6core}, with a
zoom-in around the imaginary axis on the right-hand side.

This method can be run efficiently on six cores with time steps 6.25 times larger
than those of RK$_4$.  After normalizing for number of function evaluations, the scheme
achieves time-to-solution $0.7675/0.7071 = 8.5\%$ faster than RK$_4$ but with eighth order of accuracy.  
Compared to RK$_8$, though, we achieve time-to-solution 269\% faster.
As will be seen in Sect.~\ref{sec:waveresults}, speed-up to achieve a specified accuracy
is far improved over RK$_4$ due to the size of the leading order error term combined
with eighth order convergence to the true solution.

\subsubsection{The Eighth Order, Eight Core Method: GBS$_{8,8}$}
\begin{figure}[H]
\includegraphics[width=\linewidth]{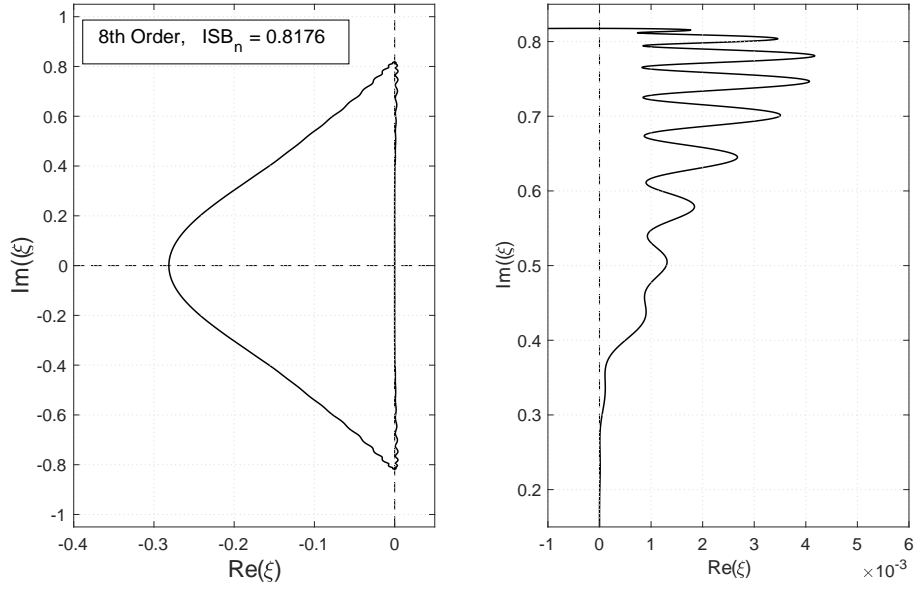}
\caption{Stability domain of the eighth order, eight core scheme (left), and zoom (right)}
\label{fig:stabdomain_p8_8core}
\end{figure}

Likewise for the six core method, we produce a robust stability domain for an eighth order
scheme that runs efficiently on eight cores, called GBS$_{8,8}$. 
The scheme's ISB$_{\text{n}}$ is 0.8176, a 0.24\% reduction from the optimal eight core value, 0.8196.
The scheme utilizes the following step counts and extrapolation weights:
\begin{equation*}
\begin{aligned}
\{n_{\text{dep}}\} &= \hspace{.3ex} \left\{\begin{matrix} 2, & 26, & 28, & 30 \end{matrix}\right\} \\
\{n_{\text{free}}\} &= \hspace{.3ex} \left\{\begin{matrix} 4, & 6, & 8, & 10, & 12, & 14, & 16, & 18, & 20, & 22 \end{matrix}\right\} \\
c_{\text{free}} &= 
\begin{matrix}[1.5]
    \bigg[ & \frac{6833}{476577792}, & \frac{10847}{91078656}, & \frac{15235}{34643968}, & \frac{383}{321152}, & \frac{543}{198784}, & \hdots
\end{matrix} \\
& \hspace{16ex}
\begin{matrix}[1.5]
    \hdots & \frac{9947}{1741056}, & \frac{6243}{543104}, & \frac{6875}{296192}, & \frac{1401}{28496}, & \frac{17713}{152688}, & \frac{6375}{19264} & \bigg]^T.
\end{matrix}
\end{aligned}
\end{equation*}
The scheme's stability domain is plotted in Fig.~\ref{fig:stabdomain_p8_8core}, with a
zoom-in around the imaginary axis on the right-hand side.

This method can be run efficiently on eight cores with time steps 8.96 times larger
than those of RK$_4$.  After normalizing for number of function evaluations the scheme
achieves time-to-solution 15.6\% faster than RK$_4$, and 287\% faster
than RK$_8$.  The two additional cores grant us a 6.5\% increase in efficiency over GBS$_{8,6}$.

\subsubsection{The Twelfth Order, Eight Core Method: GBS$_{12,8}$}
The following twelfth order method runs on eight cores and is therefore dubbed GBS$_{12,8}$.
The scheme's ISB$_{\text{n}}$ is 0.7116, a 0.17\% reduction from the optimal eight core value of 0.7128.
The scheme can be implemented using the following step counts and extrapolation weights:
\begin{equation*}
\begin{aligned}
\{n_{\text{dep}}\}  &= \hspace{.3ex} \left\{\begin{matrix} 2, & 8, & 10, & 16, & 24, & 26 \end{matrix}\right\} \\
\{n_{\text{free}}\} &= \hspace{.3ex} \left\{\begin{matrix} 4, & 6, & 12, & 14, & 18, & 20, & 22, & 28, & 30 \end{matrix}\right\} \\
c_{\text{free}} &= 
\begin{matrix}[1.5]
    \bigg[ & \frac{235}{21030240256}, & \frac{4147}{1612709888}, & \frac{11521}{39731200}, & \frac{2375}{3528704}, & \frac{6435}{708736}, & \hdots
\end{matrix} \\
& \hspace{27ex}
\begin{matrix}[1.5]
    \hdots & \frac{1291}{15780}, & \frac{11311}{4672}, & -\frac{180864}{751}, & \frac{222080}{2079} &\bigg]^T.
\end{matrix}
\end{aligned}
\end{equation*}
Stable time steps with GBS$_{12,8}$ are 7.79 times larger than those of RK$_4$
and, after normalization, time-to-solution is improved by 0.6\%.
This (very) modest efficiency improvement is drastically offset by the twelfth order of convergence
of the method - wall time to achieve a desired accuracy is far shorter than that of RK$_4$.

\begin{figure}[H]
\includegraphics[width=\linewidth]{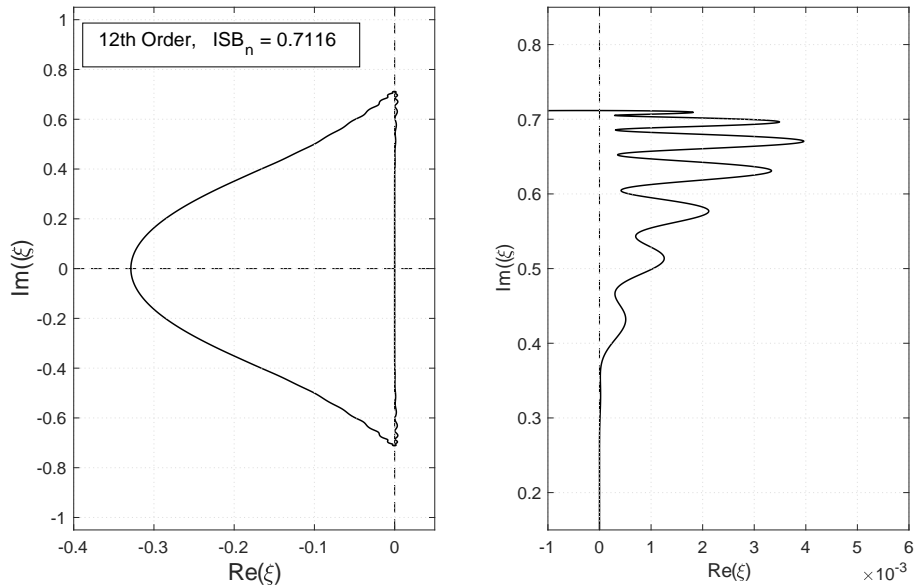}
\caption{Stability domain of the twelfth order, eight core scheme (left), and zoom (right)}
\label{fig:stabdomain_p12_8core}
\end{figure}

\section{Numerical Results}\label{sec:numerical_results}

The following results demonstrate the optimized time
steppers on a test problem with known analytic solutions.

\subsection{One-Way Wave Equation}\label{sec:waveresults}

To demonstrate the performance benefit over standard time steppers from the
literature we run the periodic one-way wave equation, 
\begin{equation}
\begin{aligned}
& \tfrac{\partial}{\partial t}u(x,t) + \tfrac{\partial}{\partial x}u(x,t) = 0, \hspace{4ex}&& 0 \le x \le 1, \hspace{2ex}t>0, \\
& u(x,0) = \tfrac{1}{2}\left(1-\cos{2\pi x}\right), && 0 \le x \le 1,
\end{aligned}
\end{equation}
utilizing the rational-coefficient GBS$_{8,6}$ and GBS$_{12,8}$ methods, with RK$_4$ as reference.
Spatial derivatives are spectral to ensure errors are due to the time stepping algorithm alone.  We run
all time steppers near their respective limits of stability, at $\lambda=\Delta t/\Delta x=\frac{.99}{\pi}\times \text{ISB}$,
where the factor of $\pi$ arises from the spectral spatial derivatives.
After convecting the wave once around the periodic interval we compute the absolute error with
respect to the analytic solution, then refine in both time and space.  

\begin{figure}[H]
\includegraphics[width=\linewidth]{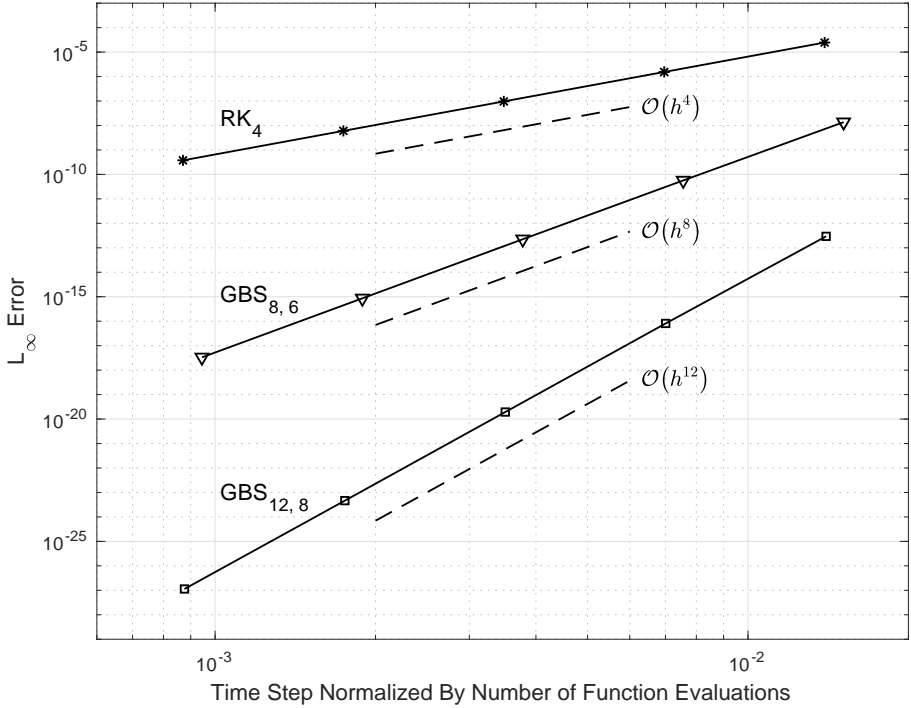}
\caption{Convection error vs. normalized time step for various methods}
\label{fig:convection_error}
\end{figure}

Convergence to the analytic solution for the various methods
is demonstrated in Fig.~\ref{fig:convection_error}.  For fair comparison across methods,
the horizontal axis is time step normalized by the number of function evaluations per step.
Thus vertical slices correspond to equal time-to-solution, neglecting the overhead of sharing
data across cores.  We use a high precision floating point library \cite{MCT15} for computation since
machine precision is achieved in the high order methods before we can establish a 
trendline.  Coefficient truncation to double precision causes error to
stagnate at $10^{-14}$ and $10^{-12}$ for the eighth and twelfth order methods, respectively.  
To obtain full floating point precision to $10^{-16}$ the extrapolation coefficients
must be precise to twenty significant digits.

\section{Conclusions}

We have presented a scheme for maximizing the time step size for extrapolation
based ODE solvers.  To do so we construct an underdetermined Vandermonde system,
then optimize the weights to maximize the stability domain along a given curve in
the complex plane.  For wave-type PDEs we utilize GBS integrators and optimize the
methods for imaginary axis coverage.  We achieve large ISB
values for methods through order sixteen which, when implemented on a computer with
several cores, yield faster time to solution than standard Runge-Kutta integrators.

The optimization method leaves both the time integrator and desired contour of stability
as user parameters.  Changing the ODE integrator in turn changes the stability polynomial
basis, which immediately affects the resulting extrapolated stability domains.  The GBS
integrator maintains large ISB through extrapolation; other integrators may be better suited
for different desired stability domains.  Future work therefore involves identifying other
suitable integrators for stability domain optimization in different contexts.  

The underdetermined extrapolation scheme saturates in parallelism
around ten cores.  We can improve scalability by incorporating the optimized schemes as local
building blocks in time-parallel solvers like Parareal \cite{LMT01}.  These solvers are known to be
less efficient with wave-type PDEs \cite{G15}.  Stable algorithms may be achieved by optimizing the
global time-parallel integrators rather than optimizing the coarse and fine grid propagators individually.
These degrees of freedom provide more flexibility than optimizing only the local schemes and is
a promising research direction for improving time to solution for wave-type equations.

\begin{appendices}
\renewcommand\thesection{Appendix \Alph{section}}
\section{Stability domains and imaginary axis coverage for some standard classes of ODE solvers}
\label{sec:appendix_a}
\renewcommand\thesection{\Alph{section}}

\subsection{Stability domains and their significance for MOL time stepping}

Each numerical ODE integration technique has an associated \emph{stability
domain}, defined as the region in a complex $\xi$-plane, with $\xi=h \lambda$,
for which the ODE method does not have any growing solutions when
it is applied to the constant coefficient ODE 
\begin{equation}
y^{\prime}=\lambda y\:.
\label{eq:SD_ODE}
\end{equation}
For a one-step method the \emph{stability polynomial}, here denoted
$R(\xi)$, is the numerical solution after one step for Dahlquist's
test equation (\ref{eq:SD_ODE}) \cite{NW91}.  The method's stability
domain is then 
\begin{align}
S=\{\xi \in \mathbb{C} : |R(\xi)| \le 1\}.
\end{align}
When solving ODEs, the stability domain can provide a guide to the
largest time step $h$ that is possible without a decaying solution
being misrepresented as a growing one. In the context of MOL-based
approximations to a PDE of the form $\frac{\partial u}{\partial t}=L(x,t,u)$,
the role of the stability domain becomes quite different, providing
necessary and sufficient conditions for numerical stability under
spatial and temporal refinement: all eigenvalues to the discretization
of the PDE's spatial operator $L$ must fall within the solver's stability
domain. For wave-type PDEs, the eigenvalues of $L$ will predominantly
fall up and down the imaginary axis. As long as the time step $h$ is
small enough, this condition can be met for solvers that feature a
positive ISB, but never for solvers with ISB = 0. 

\subsection{Runge-Kutta methods}

\begin{figure}[H]
\includegraphics[width=\linewidth]{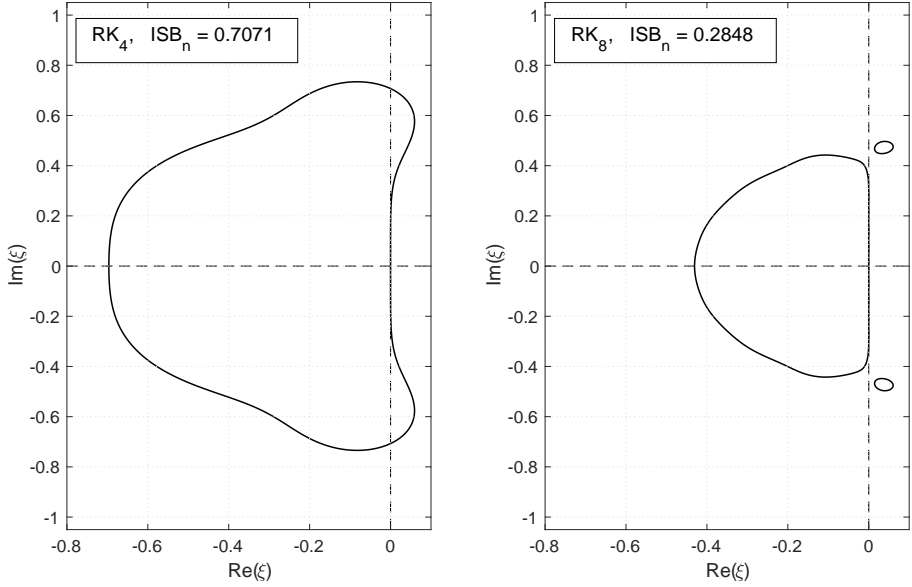}
\caption{Normalized stability domains for RK$_4$ (left) and RK$_8$ (right)}
\label{fig:stabdomain_rk}
\end{figure}

All $p$-stage RK methods of order $p$ feature the same stability
domains when $p=1,2,3,4.$ For higher orders of accuracy, more than $p$ stages
(function evaluations) are required to obtain order $p$. The $\textrm{R}\textrm{K}_{4}$
scheme used here is the classical one, and the $\textrm{R}\textrm{K}_{8}$
scheme is the one with 13 stages, introduced by Prince and Dormand
\cite{PD81}, also given in \cite{HNW87} Table 6.4. Their normalized stability
domains are shown in Fig.~\ref{fig:stabdomain_rk}.
Their ISBs are 2.8284 and 3.7023, respectively. 

\renewcommand\thesection{Appendix \Alph{section}}
\section{ISB Optimization Algorithm}
\label{sec:appendix_b}
\renewcommand\thesection{\Alph{section}}

\subsection{Optimization Formulation}
Let the extrapolated GBS stability polynomial be denoted $R(\xi)$, and the individual
stability polynomials from each of $m$ extrapolation components be denoted $P_i(\xi)$.
Then we have
\begin{align}
R(\xi) = \sum_{i=1}^m{c_iP_i(\xi)}
\end{align}
for extrapolation weights $c_i$, $1 \le i \le m$.
We collect the monomial coefficients of each $P_i(\xi)$ into the rows of a matrix, denoted $P(\xi)$.
Then $R(\xi)$ can be more compactly expressed as $R(\xi)=c^TP(\xi)$.  Now let the
left-hand-side Vandermonde matrix from (\ref{eq:vander}) be denoted $V$, and the right-hand-side
constraint vector be denoted $b$.  Then our order constraint equation (\ref{eq:vander}) can be
rewritten as $Vc=b$.

Denoting the time step size $h$, and given a curve  $\Lambda \subset \mathbb{C}$,
we specify the optimization problem as follows:
\begin{equation}
\begin{aligned}
& \underset{c_1,c_2,\ldots,c_m}{\text{maximize}}
& & h \\
& \text{subject to}
& & |R(h\lambda)|-1 \leq 0, \; \forall \lambda \in \Lambda. \\
& & & Vc = b
\end{aligned}
\label{eq:problem1}
\end{equation}

Following the work of Ketcheson and Ahmadia in \cite{KA12}, we reformulate the
optimization problem in terms of an iteration over a convex subproblem.  Minimizing
the maximum value of $|R(h\lambda)|-1$ over the weights $c_i$ is a convex problem 
(see \cite{KA12}).  We therefore define the subproblem as follows:
\begin{equation}
\begin{aligned}
& \underset{c_1,c_2,\ldots,c_m}{\text{minimize}}
& & \max_{\lambda\in\Lambda}{(|R(h\lambda)|-1)} \\
& \text{subject to}
& &  Vc = b
\end{aligned}
\label{eq:cvx_subproblem}
\end{equation}
Calling the minimax solution to (\ref{eq:cvx_subproblem}) $r(h,\Lambda)$, we can now
reformulate the optimization problem as:
\begin{equation}
\begin{aligned}
& \underset{c_1,c_2,\ldots,c_m}{\text{maximize}}
& & h \\
& \text{subject to}
& & r(h,\Lambda) \leq 0 \\
\end{aligned}
\label{eq:reformulated_problem}
\end{equation}
The optimization routine was implemented with the CVX toolbox for MATLAB \cite{CVX14}
using a bisection over time step $h$.  Results presented in this paper use the software
\emph{OPTISB} \cite{AE19} to optimize the stability domains.

\subsection{Comparison to Optimizing Monomial Coefficients}

The main theoretical difference between our current algorithm and the algorithm
presented in \cite{KA12} is the basis over which coefficients are
optimized.  In \cite{KA12} the authors optimize directly the coefficients
to the stability polynomial in the monomial basis.
This yields an optimal stability polynomial that must be approximated
with a Runge-Kutta integrator.  The polynomial is therefore fed into a second
optimization routine to compute the Runge-Kutta coefficients.

In the extrapolation coefficient optimization we operate directly on linear
combinations of the time stepper stability polynomials.  The true optimal stability
polynomial therefore may not be in the space of extrapolated GBS time stepper
stability polynomials.  However, the resulting stability polynomial is immediately
realizable and we require no further optimization stage to generate our time
stepping algorithm.  

\subsection{Implementing Order Constraints}

To guarantee accuracy we require order constraints to be satisfied to machine
precision.  Most optimization routines accept equality constraints that will hold within a
certain tolerance.  Due to ill-conditioning of the Vandermonde systems we
prefer to explicitly enforce the order constraints in the convex optimization.
As in Sect.~\ref{sec:notation} we split the stability polynomials into two groups which take on the ``dep" and
``free" subscripts, denoting dependent and optimized quantities, respectively.  The dependent
weights guarantee the extrapolation scheme achieves the specified order of accuracy.
The remaining weights are our optimization variables.   Thus the stability
polynomial is computed as follows:
\begin{align}
R(\xi) = c_{\text{dep}}^TP_{\text{dep}}(\xi) + c_{\text{free}}^TP_{\text{free}}(\xi).
\end{align}
Order constraints take the form (\ref{eq:order_constraints}) which yields the
dependent weight computation (\ref{eq:cdep_computation}).
Splitting the weights apart reduces the number of design variables and, in practice, leads to
better solutions than when utilizing equality constraints.

\end{appendices}

\bibliographystyle{spmpsci}
\bibliography{References.bib}

\end{document}